\documentclass[11pt]{amsart}
\usepackage[latin1]{inputenc}
\usepackage{a4,amssymb}
\usepackage{enumerate}

\usepackage[T1]{fontenc}
\usepackage{mathptmx}

\usepackage[colorlinks, bookmarks=true]{hyperref}

\newtheorem{theorem}{Theorem}[section]

\newtheorem{lemma}[theorem]{Lemma}
\newtheorem{corollary}[theorem]{Corollary}
\theoremstyle{definition}
\newtheorem{remark}[theorem]{Remark}

\newtheorem{problem}[theorem]{Problem}
\numberwithin{equation}{section}

\newcommand{\unit}{\mathbf{1}}

\begin{document}
	
\title[Lie--Trotter formulae in Jordan-Banach algebras]{Lie--Trotter formulae in Jordan--Banach algebras
 with applications to the study of spectral-valued multiplicative functionals}

\author[G.M. Escolano]{Gerardo M. Escolano}

\address{Departamento de An{\'a}lisis Matem{\'a}tico, Facultad de Ciencias, Universidad de Granada \hyphenation{Gra-nada}, 18071 Granada, Spain.}
\email{gemares@correo.ugr.es}

\author[A.M. Peralta]{Antonio M. Peralta}

\address{Departamento de An{\'a}lisis Matem{\'a}tico, Facultad de Ciencias, Universidad de Granada, 18071 Granada, Spain.
Instituto de Matem{\'a}ticas de la Universidad de Granada (IMAG).}
\email{aperalta@ugr.es}

\author[A.R. Villena]{Armando R. Villena}

\address{ Departamento de An{\'a}lisis Matem{\'a}tico, Facultad de Ciencias, Universidad de Granada, 18071 Granada, Spain.
Instituto de Matem{\'a}ticas de la Universidad de Granada (IMAG).}
\email{avillena@ugr.es}



\maketitle

\begin{abstract} We establish some Lie--Trotter formulae for unital complex Jordan--Banach algebras, showing that for each couple of elements $a,b$ in a unital complex Jordan--Banach algebra $\mathfrak{A}$ the identities 
$$ 	\lim_{n\to \infty} \left(e^{\frac{a}{n}}\circ e^{\frac{b}{n}} \right)^{n} = e^{a+b},\ \lim_{n\to \infty} \left(U_{e^{\frac{a}{n}}} \left( e^{\frac{b}{n}}\right)  \right)^{n} = e^{2 a+b}, \hbox{ and }$$	$$  \lim_{n\to \infty} \left(U_{e^{\frac{a}{n}},e^{\frac{c}{n}}} \left( e^{\frac{b}{n}}\right)  \right)^{n} = e^{a+b + c}$$ hold. These formulae are actually deduced from a more general result involving holomorphic functions with values in $\mathfrak{A}$. These formulae are employed in the study of spectral-valued (non-necessarily linear) functionals $f:\mathfrak{A}\to \mathbb{C}$ satisfying $f(U_x (y))=U_{f(x)}f(y),$ for all $x,y\in \mathfrak{A}$. We prove that for any such a functional $f,$ there exists a unique continuous (Jordan-)multiplicative linear functional $\psi\colon \mathfrak{A}\to\mathbb{C}$
such that $ f(x)=\psi(x),$ for every $x$ in the connected component of set of all invertible elements of $\mathfrak{A}$ containing the unit element. If we additionally assume that $\mathfrak{A}$ is a JB$^*$-algebra and $f$ is continuous, then $f$ is a linear multiplicative functional on $\mathfrak{A}$. The new conclusions are appropriate Jordan versions of results by Maouche, Brits, Mabrouk, Shulz, and Tour{\'e}.
\end{abstract}

\section{Introduction}

Current research in functional analysis, and especially on preservers, is particularly interested on exploring when the traditional conditions of linearity can be somehow removed in classic results at the cost of assuming the preservation of another algebraic or geometric property. This is the case for the celebrated  Gleason-Kahane-\.{Z}elazko theorem, asserting that a non-zero linear functional $\phi : A\to \mathbb{C},$ where $A$ is an (associative) complex Banach algebra is multiplicative if and only if it is spectral-valued, that is, $\phi (a)$ belongs to the spectrum of $a$, $\sigma (a),$ in $A$ (see \cite[Theorem 4.1.1]{AupetitBook}). The following natural problem on preservers is directly motivated by the just recalled theorem.

\begin{problem}\label{Problem reciprocal of the GKZ towards linearity} Let $f:A\to \mathbb{C}$ be a multiplicative functional from a (unital) complex Banach algebra. Suppose that for each $x\in A$, $f(x)$ lies in the spectrum of $x$. Is $f$ a multiplicative linear functional on $A$?  
\end{problem}

This problem was considered by Maouche in his pioneering work \cite{Maouche1996}, where he established the following:

\begin{theorem}\label{t Maouche}{\rm \cite{Maouche1996}}
Let $A$ be a unital complex Banach algebra, and let $f\colon A\to \mathbb{C}$ be a (non-necessarily continuous) multiplicative spectral-valued functional.  Then there exists a unique multiplicative linear functional $\psi$ on $A$ that agrees with $f$ on the connected component of the identity in the group of all invertible elements of $A$.
\end{theorem} 

The conclusion of the previous theorem is optimal in the sense that $f$ need not be, in general, linear \cite[Exemple]{Maouche1996}. After Maouche's contribution positive solutions to Problem~\ref{Problem reciprocal of the GKZ towards linearity} have been found for commutative unital C$^*$-algebras by Tour{\'e} and Brits \cite{TouBrits2020}, von Neumann algebras by Tour{\'e}, Shulz and Brits \cite{TouShulzBrits2017}, and for unital C$^*$-algebras by Brits, Mabrouk and Tour{\'e} \cite{BritsMabroukToure2021}. The principal aim of this paper is to explore Problem~\ref{Problem reciprocal of the GKZ towards linearity} in the case of unital complex Jordan--Banach algebras. However, in order to understand the appropriate form of this problem in the Jordan setting we shall review some of the basic notions and additional tools developed in this study.   

Two key results appear among the arguments developed to address the problem in the aforementioned references, namely the Lie--Trotter formula and the Kowalski-S{\l}odkowski theorem (see section~\ref{sec: GKZ and Maouche} for more details on the latter theorem). 

The so-called \emph{Lie--Trotter formula}, a celebrated result in operator theory, asserts that for any arbitrary elements $a,b$ in an associative complex Banach algebra $A$ we have \begin{equation}\label{eq Lie-Trotter formula Banach algebras}
 \lim_{n\to \infty} \left(e^{\frac{a}{n}} e^{\frac{b}{n}} \right)^{n} = e^{a+b},	
\end{equation} where we write $e^x$ for the usual exponential of an element $x$ in $A$ (see, for example, \cite[Exercise 15, page 67]{AupetitBook}, \cite[Theorem 2.10]{HallBook2003}).  The Lie--Trotter formula holds even for pairs of non-commuting elements. Applications of this result can be found in quantum theory \cite{Simon79,Davies80,HallBook2013}. For example, the Feynman--Kac formula can be obtained via the Lie--Trotter formula \cite[\S 6.3.1]{ApplebaumBook2019}.

A real or complex \emph{Jordan--Banach algebra} is a Banach space $\mathfrak{A}$ equipped with a bilinear mapping $(a,b)\mapsto a\circ b$ (called the Jordan product) satisfying the following axioms:\begin{enumerate}[$(J1)$]\item $\|a\circ b\|\leq \| a\| \, \|b\|,$  for all $a,b\in \mathfrak{A}$;
\item $a\circ b = b\circ a,$ for all $a,b\in \mathfrak{A}\ \ $ (commutativity);
\item $(a^2\circ b)\circ a = (a\circ b)\circ a^2,$ for all $a,b\in \mathfrak{A}\ \ $ (Jordan identity).
\end{enumerate}	Let us observe that associativity of the product has been replaced by the Jordan identity. In this note we shall only consider complex Jordan--Banach algebras. We shall say that $\mathfrak{A}$ is \emph{unital} if there exists an element $\mathbf{1}\in \mathfrak{A}$ (called the unit of $\mathfrak{A}$) such that $\mathbf{1}\circ a= a$ for all $a\in \mathfrak{A}$. It is hard to understand well Jordan algebras without the $U$ maps defined as follows.  Given an element $a \in \mathfrak{A}$ the symbol $U_a$ will stand for the linear mapping on $\mathfrak{A}$ defined by $U_a (b) := 2(a\circ b)\circ a - a^2\circ b$ for each $b\in\mathfrak{A}$.

A natural example can be derived from any associative Banach algebra $A$ equipped with the natural Jordan product given by $a\circ b := \frac12 ( a b + ba )$. Jordan--Banach subalgebras of Banach algebras are called \emph{special}. The hermitian part of each C$^*$-algebra is a real Jordan--Banach subalgebra but not an associative subalgebra. Although in this note the Jordan-Banach algebras are studied from a purely mathematical point of view, we observe that Jordan algebras were introduced with the aim of being a suitable setting for axiomatic quantum mechanics by P. Jordan, J. von Neumann and E. Wigner \cite{JorvNeuWign34}.

Having in mind that the closed Jordan subalgebra, $\mathfrak{A}_a$, generated by a single element $a$ (and the unit element) of a complex Jordan--Banach algebra, $\mathfrak{A},$ is an associative and commutative complex Banach algebra, the usual \emph{holomorphic functional calculus} makes sense in a very natural way in the setting of complex Jordan--Banach algebras (cf. \cite[\S2, Theorem 2.3]{Aupetit95} and  \cite[Theorems 4.1.88 and 4.1.93]{CabGarPalVol1}). An element $a$ in a unital Jordan--Banach algebra $\mathfrak{A}$ is called \emph{invertible} if we can find  $b\in \mathfrak{A}$ satisfying $a \circ b = \mathbf{1}$ and $a^2 \circ b = a.$ The element $b$ is unique and it will be denoted by $a^{-1}$ (cf. \cite[3.2.9]{HOS} or \cite[Definition 4.1.2]{CabGarPalVol1}). The symbol $\mathfrak{A}^{-1}$ will stand for the set of all invertible elements in $\mathfrak{A}$.  Having these basic notions in mind, it is clear that we can define the exponential of each element $a$ in $\mathfrak{A}$ by the holomorphic functional calculus or by the power series $\displaystyle e^{a} =\exp (a)= \sum_{n=0}^{\infty} \frac{a^n}{n!}$, which gives an invertible element in $\mathfrak{A}$ and lies in the Jordan--Banach subalgebra generated by ${a}$ and $\mathbf{1}$. A basic spectral theory for complex Jordan--Banach algebras satisfying  similar properties to those in the setting of associative complex Banach algebras is fully available (see \cite[Chapter 4]{CabGarPalVol1} or the next subsection~\ref{subsec:background} for more details). The Jordan-spectrum of an element $a\in \mathfrak{A}$ is the (non-empty compact) set $\hbox{J-}\sigma (a):=\{\lambda\in \mathbb{C}: a-\lambda\unit\notin\mathfrak{A}^{-1}\}.$ 

The Lie--Trotter formula remains completely unexplored in the setting of complex Jordan--Banach algebras. We should begin by asking ``What is the most appropriate Jordan version of the Lie--Trotter formula in \eqref{eq Lie-Trotter formula Banach algebras}?'' The first candidate comes by replacing the associative product by the Jordan one, that is, 
\begin{equation}\label{eq Lie-Trotter formula Jordan product}
	\lim_{n\to \infty} \left(e^{\frac{a}{n}}\circ e^{\frac{b}{n}} \right)^{n} = e^{a+b},	\ \ \ \ (a,b\in \mathfrak{A}).
\end{equation} 
This version seems interesting and newfangled even in the associative setting. However, having in mind that the Jordan product of two invertible elements is not, in general, an invertible element -- actually there exists elements of the form $e^{a}$ and $e^b$ whose Jordan product is zero-- but expressions of the form $U_a (b) := 2(a\circ b)\circ a - a^2\circ b$ give invertible elements whenever $a$ and $b$ are, the second natural formula might be following:
\begin{equation}\label{eq Lie-Trotter formula Jordan U product}
	\lim_{n\to \infty} \left(U_{e^{\frac{a}{n}}} \left( e^{\frac{b}{n}}\right)  \right)^{n} = e^{2 a+b}, \hbox{ or more generally, } \lim_{n\to \infty} \left(U_{e^{\frac{a}{n}},e^{\frac{c}{n}}} \left( e^{\frac{b}{n}}\right)  \right)^{n} = e^{a+b + c},
\end{equation} for all $a,b, c\in \mathfrak{A},$ where $U_{a,c} (b) := (a\circ b) \circ c + (c\circ b) \circ a - (a\circ c) \circ b$ with $U_{a,a} (b) = U_{a} (b).$ We prove in section~\ref{sec: Lie--Trotter formulae} that all the formulae in \eqref{eq Lie-Trotter formula Jordan product} and \eqref{eq Lie-Trotter formula Jordan U product} hold true in the setting of unital complex Jordan--Banach algebras. We actually prove a much general result showing that if $f\colon \mathbb{D}_r\to \mathfrak{A}$ is a holomorphic mapping from the open disc of radius $r$ in $\mathbb{C}$ into a unital complex Jordan--Banach algebra with $f(0)=\mathbf{1}$, $(\lambda_n)$ is a sequence in $\mathbb{D}_r$ with $\lim_{n\to\infty}\lambda_n=0$, and $(\mu_n)$ is a sequence in $\mathbb{C}$  such that $\lim_{n\to\infty} \lambda_n\mu_n=\lambda\in\mathbb{C}$, then for each sufficiently large $n\in\mathbb{N}$, $\bigl(f(\lambda_n)\bigr)^{\mu_n}$ is well-defined
	and
	$$
	\lim_{n\to\infty}\bigl(f(\lambda_n)\bigr)^{\mu_n}=
	\exp\bigl(\lambda f'(0)\bigr)
	$$ (see Theorem~\ref{t Lie-Trotter fla with holomorpic mappings}). The validity of formulae \eqref{eq Lie-Trotter formula Jordan U product} and \eqref{eq Lie-Trotter formula Jordan product} is established as a consequence in Corollary~\ref{c Li-Trotter formulae Jordan-Banach algeras}.
	
The beginning of section~\ref{sec: GKZ and Maouche} is devoted to state a version of the celebrated Gleason-Kahane-\.{Z}elazko for complex Jordan--Banach algebras which is deduced from the corresponding version for associative complex Banach algebras. More concretely, a non-zero linear functional on a complex Jordan--Banach algebra $\mathfrak{A}$ is (Jordan-)multiplicative if and only if $\varphi (a)\in \hbox{J-}\sigma(a)$ for all $a\in \mathfrak{A}$ (see Theorem~\ref{t GKZ JB}). Assuming next that $\mathfrak{A}$ is unital, our next goal will consist in exploring how close is a spectral-valued Jordan--multiplicative functional $f: \mathfrak{A}\to \mathbb{C}$ from agreeing with a Jordan-multiplicative functional on the principal component of the invertible elements, in the line of the result proved by Maouche for Banach algebras \cite{Maouche1996}. A milestone result in the theory of Jordan--Banach algebras, due to Loos \cite{Loos96}, provides an algebraic characterisation of the principal component, $\mathfrak{A}_{\mathbf{1}}^{-1}$ (i.e., the connected component containing $\unit$) of the invertible elements, $\mathfrak{A}^{-1},$ in a unital complex Jordan--Banach, $\mathfrak{A},$ in the form $$\mathfrak{A}_{\mathbf{1}}^{-1} = \left\{ U_{\exp(a_1)}\cdots U_{\exp(a_n)} (\textbf{1}) : a_i\in \mathfrak{A}, \ n\in \mathbb{N} \right\}.$$ The special form of this principal connected component of $\mathfrak{A}^{-1},$ and the lacking of an argument proving that a (non-necessarily linear) Jordan-multiplicative mapping preserves triple products of the form $U_{a} (b)$, forces us to consider different hypotheses in the following Jordan version of Problem~\ref{Problem reciprocal of the GKZ towards linearity}. 

\begin{problem}\label{Problem reciprocal of the GKZ towards linearity Jordan-Banach} Let $f:\mathfrak{A}\to \mathbb{C}$ be a (non-necessarily linear) functional satisfying:
\begin{enumerate}[$(a)$]
		\item $f(U_x (y))=U_{f(x)}f(y),$ for all $x,y\in \mathfrak{A}$,
		\item $f(x)\in\text{\rm J-}\sigma (x),$  for all $x\in \mathfrak{A}.$
\end{enumerate} Does $f$ coincide with a (Jordan-)multiplicative linear functional on the principal component of $\mathfrak{A}^{-1}$? If we also assume that $f$ is continuous, is $f$ a (Jordan-)multiplicative linear functional on $\mathfrak{A}$?  
\end{problem}

In Theorem~\ref{t Maouche Jordan Banach with natural hypotheses} we give a complete positive answer to the first question in Problem~\ref{Problem reciprocal of the GKZ towards linearity Jordan-Banach}.  Compared with the result established by Manouche in \cite{Maouche1996} for unital complex Banach algebras, our arguments do not make any use of the Kowlaski--S{\l}odkowski theorem but replace it with a weaker technical lemma. We shall also study spectral-valued maps $f\colon \mathfrak{A}\to\mathbb{C}$ which are Jordan-multiplicative (see Theorem~\ref{t Maounche Jordan mutiplicative}), and show that Maouche's result can be derived from the Jordan version thanks to a classical contribution by \.{Z}elazko around the Gleason--Kahane--\.{Z}elazko theorem. 

All the previous results on preservers are in some sense optimal, since there exists a counter-example by Maouche showing the existence of a (non-continuous) multiplicative spectral-valued functional $f : C[0,1]\to \mathbb{C}$ which is not linear (see \cite[Exemple]{Maouche1996}).

In the final section~\ref{sec: Brits} we provide a complete positive solution to the second question in Problem~\ref{Problem reciprocal of the GKZ towards linearity Jordan-Banach} in the case of unital JB$^*$-algebras. As well as C$^*$-algebras are those complex Banach algebras $A$ equipped with an algebra involution $a\mapsto a^*$ satisfying the ``geometric'' Gelfand--Naimark axiom $\| a^* a\| = \|a\|^2\ \ $ ($a\in A$), a JB$^*$-algebra $\mathfrak{A}$ is a complex Jordan--Banach algebra equipped with an algebra involution $a\mapsto a^*$ satisfying the ``geometric'' axiom $\left\| U_{a} (a^*) \right\| = \|a\|^3$ ($a\in \mathfrak{A}$). It is worth to note that if a C$^*$-algebra $A$ is equipped with its natural structure of Jordan algebra, we have $U_a (a^* ) = a a^* a$, so the geometric axiom in the definition of JB$^*$-algebras is an equivalent version of the Gelfand--Naimark axiom.  Our advances extend to the Jordan setting the most general known conclusions for C$^*$-algebras obtained by Tour{\'e} and Brits in \cite{TouBrits2020}, namely, every continuous and spectral-valued functional $f$ from a unital JB$^*$-algebra $\mathfrak{A}$ into $\mathbb{C}$ satisfying $f(U_x (y))=U_{f(x)}f(y),$ for all $x,y\in \mathfrak{A}$ is a linear (Jordan-)multiplicative functional on $\mathfrak{A}$ (see Theorem~\ref{t BMT for JB}). Our approach is independent from the arguments in the setting of C$^*$-algebras, and it can be said that a slightly simpler.

\subsection{Background on Jordan--Banach algebras}\label{subsec:background}

\emph{Power associativity} is one of the basic properties of Jordan algebras, more concretely, for each Jordan-Banach algebra $\mathfrak{A}$, the Jordan--Banach subalgebra, $\mathfrak{A}_{a},$ generated by a single element $a$ in a Jordan--Banach algebra, $\mathfrak{A}$, is associative, equivalently, by setting $a^0 = \textbf{1}$, and $a^{n+1} = 
a\circ a^n,$ we have $a^n \circ a^m = a^{n+m}$ for all $n,m\in \mathbb{N}\cup \{0\}$ (cf. \cite[Lemma 2.4.5]{HOS} or \cite[Corollary 1.4]{AlfsenShultz2003}). Therefore $\mathfrak{A}_{a}$ is a commutative and associative Banach algebra. 

Given elements $a,c$ in a Jordan--Banach algebra $\mathfrak{A},$ the symbol $U_{a,c}$ will stand for the linear mapping on $\mathfrak{A}$ defined by $U_{a,c} (b) := (a\circ b)\circ c + (c\circ b)\circ a - (a\circ c)\circ b.$ Clearly $U_{a,a}$ agrees with the mapping $U_a$ mentioned in the introduction. Invertibility is deeply related to $U$ maps since an element $a\in \mathfrak{A}$ is invertible if and only if $U_a$ is a bijective mapping, and in such a case $U_a^{-1} = U_{a^{-1}}$  (cf. \cite[3.2.9]{HOS} and \cite[Definition 4.1.2 and Theorem 4.1.3]{CabGarPalVol1}).
We shall write $\mathfrak{A}^{-1}$ for the set of all invertible elements in $\mathfrak{A}$. It is known that $\mathfrak{A}^{-1}$ is open (see \cite[Theorem 4.1.7]{CabGarPalVol1}), and hence $\mathfrak{A}^{-1}$ is locally connected and all its connected components are open. 

The Jordan product of two invertible elements is not, in general, an invertible element. Anyway, an element of the form $U_x (y)$ is invertible if and only if both $x$ and $y$ are (cf. \cite[Theorem 4.1.3$(vi)$]{CabGarPalVol1}), and in such a case $(U_x (y))^{-1} = U_x^{-1} (y^{-1})$; hence, $(x^{-1})^2 = (x^2)^{-1} = x^{-2}$. Consequently, the square and the $n$-th power ($n\in \mathbb{N}$) of each invertible element is an invertible element. 

The \emph{spectrum} (also called the \emph{Jordan-spectrum}) of an element $a$ in a unital complex Jordan-Banach algebra $\mathfrak{A}$ (denoted by J-$\sigma(a)$) is the set of all $\lambda\in \mathbb{C}$ for which $a-\lambda \textbf{1}$ is not invertible. As in the setting of associative Banach algebras, J-$\sigma(a)$ is a non-empty compact subset of the complex plane (see \cite[Theorem 4.1.17]{CabGarPalVol1}. We shall frequently apply along this paper that if we denote by $A^{sym}$ the Jordan--Banach algebra naturally associated to a unital complex Banach algebra $A$, an element $a$ in $A$ is invertible in the Jordan sense if and only if it is invertible in the usual sense.  In particular, $\sigma (a) = \hbox{J-}\sigma (a)$ for every $a\in A.$

We refer to the monographs \cite{HOS, AlfsenShultz2003} and \cite{CabGarPalVol1} as authoritative sources on Jordan--Banach algebras and JB$^*$-algebras. For additional details on the holomorphic functional calculus the reader can consult \cite[\S2, Theorem 2.3]{Aupetit95} and \cite[Theorems 4.1.88 and 4.1.93]{CabGarPalVol1}.

\section{Lie--Trotter formula in Jordan--Banach algebras}\label{sec: Lie--Trotter formulae}

This section is the technical core of the paper and contains several Lie--Trotter type formulae in the setting of Jordan-Banach algebras. It is perhaps worth to begin our approach by showing that in the case of associative unital complex Banach algebras our Lie--Trotter formula for the Jordan product in \eqref{eq Lie-Trotter formula Jordan U product} can be deduced from the classical formula. Namely, let $A$ be a unital complex Banach algebra. It is not hard to see by induction that the identity $$ \left( e^{\frac{a}{n}} e^{\frac{b}{n}} e^{\frac{a}{n}}\right)^n =e^{-\frac{a}{n}} \left(e^{\frac{2 a}{n}}e^{\frac{b}{n}}\right)^n e^{\frac{a}{n}} $$ holds for all $a,b\in A$ and all natural $n$. So, by taking limit in the above expression and applying the original Lie--Trotter formula and the joint continuity of the product we get \begin{equation}\label{eq Lie--Trotter U special Jordan} \lim_{n\to \infty} \left( U_{e^{\frac{a}{n}}} \left(e^{\frac{b}{n}} \right)\right)^n  = \lim_{n\to \infty} \left( e^{\frac{a}{n}} e^{\frac{b}{n}} e^{\frac{a}{n}}\right)^n =\lim_{n\to \infty} e^{-\frac{a}{n}} \left(e^{\frac{2 a}{n}}e^{\frac{b}{n}}\right)^n e^{\frac{a}{n}}=e^{2 a + b}.
\end{equation} The reader might expect that if we fix two arbitrary elements $a,b$ in a unital complex Jordan--Banach algebra $\mathfrak{A}$, by the Shirshov--Cohn theorem \cite[Theorem 2.4.14]{HOS}, the Jordan subalgebra $\mathfrak{B}$ of $\mathfrak{A}$ generated by $a,b$ and $\unit$ (in a purely algebraic sense) is special, and hence we can find an associative algebra $A$ satisfying that  $\mathfrak{B}$ is a Jordan subalgebra of $A$. In principle we cannot apply the conclusion in \eqref{eq Lie--Trotter U special Jordan} since the Shirshov--Cohn theorem is a purely algebra result and we do not know, in general, whether the associative algebra $A$ can be assumed to be a Banach algebra. Further, even in the associative case, it does not seem that \eqref{eq Lie-Trotter formula Jordan product} can be obtained from the classical Lie-Trotter formula.

We shall next employ holomorphic vector-valued functions in the standard sense. Let $X$ be a complex Banach space, a function $f$ from a nonempty open set $\Omega\subseteq \mathbb{C}$ is called (\emph{strongly}) \emph{complex-differentiable} or \emph{holomorphic} if $\displaystyle \lim_{\lambda\to \lambda_0}
\frac{f (\lambda) - f (\lambda_0)}{\lambda-\lambda_0}$ exists (in $X$) for all $\lambda_0\in \Omega$. It is known that a mapping $f: \Omega \to X$ is holomorphic if and only if it is \emph{weakly holomorphic}, that is,  for each functional $\varphi\in X^*$ the function $\varphi \circ f: \Omega\to \mathbb{C}$ is holomorphic \cite[Corollary 15.2.2]{GarretBook2019}. Each holomorphic function $f$ satisfies the Cauchy integral formula given in terms of the a Gelfand--Pettis $X$-valued
integral \cite[Corollary 15.2.2]{GarretBook2019}, actually $f$ satisfies the usual Cauchy-theory integral formulae and it is, in particular, infinitely differentiable and expressible as a convergent power series with coefficients given by Cauchy's
formulae \cite[Corollary 15.2.3]{GarretBook2019}. We shall also deal with holomorphic maps from an open subset $\Omega$ of a complex Banach space $X$ with values in another complex Banach space $Y$ as those mappings $f:\Omega \to Y$ which are Fréchet differentiable at every point in $\Omega$. Examples of the latter can be given by the holomorphic functional calculus in complex Banach and Jordan--Banach algebras (cf. \cite[Theorems 1.3.21 and 4.1.93 ]{CabGarPalVol1}).

We can now state the main result of this section. Henceforth we shall write $\mathbb{D}_r$ for the open disc in $\mathbb{C}$ of centre $0$ and radius $r>0$. The symbol $\mathbb{D}$ will stand for  $\mathbb{D}_1$. 

\begin{theorem}\label{t Lie-Trotter fla with holomorpic mappings}
Let $\mathfrak{A}$ be a unital complex  Jordan-Banach algebra and let $f\colon \mathbb{D}_r\to \mathfrak{A}$ be 
a holomorphic mapping such that $f(0)=\mathbf{1}$.	Suppose $(\lambda_n)$ is a sequence in $\mathbb{D}_r$ with 
$\lim_{n\to\infty}\lambda_n=0$, and $(\mu_n)$ is a sequence in $\mathbb{C}$  such that
 $\lim_{n\to\infty} \lambda_n\mu_n=\lambda\in\mathbb{C}$. 	
Then for each sufficiently large $n\in\mathbb{N}$, $\bigl(f(\lambda_n)\bigr)^{\mu_n}$ is well-defined
	and
	$$
	\lim_{n\to\infty}\bigl(f(\lambda_n)\bigr)^{\mu_n}=
	\exp\bigl(\lambda f'(0)\bigr).
	$$
\end{theorem}

\begin{proof} By the basic properties of the holomorphic functional calculus, the continuity of $f$ and the other hypotheses $\bigl[f(\lambda_n)\bigr]^{\mu_n}$ is well defined for $n$ large enough and we have $$
	\bigl[f(\lambda_n)\bigr]^{\mu_n}=
	\exp\bigl[\mu_n\log\,f(\lambda_n)\bigr]=
	\exp\bigl(\lambda_n\mu_n g(\lambda_n)\bigr),$$
	where $g\colon\mathbb{D}_r\to J$ is defined by
	$$
	g(\zeta)=
	\left\{
	\begin{array}{cl}
		\displaystyle{\frac{\log\bigl(f(\zeta)\bigr)}{\zeta}},&\text{if }\zeta\ne 0,\\
		f'(0),&\text{if }\zeta=0.
	\end{array}
	\right.$$
	
It suffices to show that
	$$
	\lim_{\zeta\to 0}g(\zeta)=g(0)= f'(0).$$
	By \cite[Theorem 4.1.93]{CabGarPalVol1}, the set $\Omega=\{x\in \mathfrak{A} \colon \text{J-}\sigma(x)\subset\mathbb{C}\setminus\mathbb{R}_0^-\}$
	is open and the mapping
	\[
	\Theta\colon\Omega\to
	\mathfrak{A},\quad\Theta(x)=\log(x)\quad\forall x\in\Omega
	\] 
	is holomorphic. Furthermore, given $0<\varrho<1$, the Fr\'echet differential $D\Theta(\mathbf{1})\colon \mathfrak{A}\to \mathfrak{A}$ of $\Theta$ 
	at $\mathbf{1}$ is given by
	\begin{equation*}
		\begin{split}
			D\Theta(\mathbf{1})(x)
			&=
			\frac{1}{2\pi\,i}\int_{\vert \zeta-1\vert=\varrho}\log(\zeta)\,U_{(\zeta\mathbf{1}-\mathbf{1})^{-1}}(x)\,d\zeta\\
			&=
			\frac{1}{2\pi\,i}\int_{\vert\zeta-1\vert=\varrho}\frac{\log(\zeta)}{(\zeta-1)^2}\,x\,d\zeta\\
			&=\log'(1)x=x,
			\hbox{ for all } x\in \mathfrak{A}.
		\end{split}
	\end{equation*}
	Then $U=f^{-1}(\Omega)$ is an open neighbourhood of $0$ in $\mathbb{C}$ and the mapping $\Theta\circ f$ 
	is holomorphic on $U$ with
	$$
		\lim_{\zeta\to 0}g(\zeta)= (\Theta\circ f)'(0)=D\Theta(f(0))(f'(0))=D\Theta(\unit)(f'(0))= f'(0),
	$$ as desired.
\end{proof}

The desired Lie--Trotter formulae predicted in \eqref{eq Lie-Trotter formula Jordan product} and \eqref{eq Lie-Trotter formula Jordan U product} can be now obtained as a consequence of our previous theorem. 

\begin{corollary}\label{c Li-Trotter formulae Jordan-Banach algeras}
Let $\mathfrak{A}$ be a unital complex  Jordan--Banach algebra. Then the formulae 
	$$ \lim_{n\to \infty} \left(U_{e^{\frac{a}{n}},e^{\frac{c}{n}}} \left( e^{\frac{b}{n}}\right)  \right)^{n} = e^{a+b + c},\ \ \lim_{n\to \infty} \left(U_{e^{\frac{a}{n}}} \left( e^{\frac{b}{n}}\right)  \right)^{n} = e^{2 a+b},
	$$
	$$ \hbox{ and }	\lim_{n\to \infty} \left( {e^{\frac{a}{n}}} \circ   e^{\frac{b}{n}}  \right)^{n} = e^{ a+b}, $$ hold for all $a,b,c\in \mathfrak{A}$.
\end{corollary}

\begin{proof} To prove the first formula it suffices to apply Theorem~\ref{t Lie-Trotter fla with holomorpic mappings} with
$$ f(\zeta)=U_{e^{\zeta\,  {a}},e^{\zeta\,  {c}}} \left( e^{\zeta\,  {b}}\right)
,\quad\forall \zeta\in\mathbb{C},$$ $(\lambda_n) = (\frac1n)$ and $(\mu_n) =(n)$. The second and third formulae follow from the first one by taking $a=c$ and $c=0$, respectively.
\end{proof}

\begin{remark}\label{remark further formulae} Theorem~\ref{t Lie-Trotter fla with holomorpic mappings} is a formidable device 
to produce new Lie--Trotter type formulae since there is a wide freedom on the holomorphic mapping $f$. For example, by taking 
$f(\zeta) := e^{\zeta\, a} \circ (\unit + \zeta\, b + \zeta^3\, d^3)\circ \cos (\zeta\, c),$ with $a,b,c,d$ in a Jordan--Banach 
algebra $\mathfrak{A}$ we get 
$$\lim_{n\to \infty} \left( e^{\lambda_n\, a} \circ (\unit + \lambda_n\, b + \lambda_n^3\, d^3)\circ \cos (\lambda_n\, c) \right)^{\mu_n}= 
e^{\lambda (a +b)},$$ 
for each null sequence $(\lambda_n)$ and each scalar sequence $(\mu_n)$ with $ (\lambda_n\mu_n )\to \lambda$, where the terms in the left-hand-side sequence of the formula are considered for $n$ large enough.
\end{remark}

It is worth to note that Theorem~\ref{t Lie-Trotter fla with holomorpic mappings} is also valid when $A$ is an associative unital complex Banach algebra. In such a case the holomorphic mapping $f$ can be also defined in terms of the associative product. 	

\section{Gleason--Kahane--\.{Z}elazko theorem and its reciprocal in Jordan--Banach algebras}\label{sec: GKZ and Maouche}

The aim of this section is to provide an application of the Lie--Trotter formulae established in the previous section. As we commented in the introduction, the Lie--Trotter formula for unital complex Banach algebras is employed in the study of multiplicative spectral-valued maps from a unital complex Banach algebra $A$ to the complex field by Maouche in \cite{Maouche1996}. The result affirms that for each (non-necessarily linear) multiplicative functional $\Delta: A\to \mathbb{C}$ satisfying $\Delta (a)\in \sigma (a)$, there exists a multiplicative linear functional $\psi$ on $A$ which agrees with $\Delta$ on the principal connected component of the invertible elements. Our aim here is to explore a similar conclusion for spectral-valued multiplicative functionals from a unital complex Jordan--Banach algebra into $\mathbb{C}$. 

In Maouche's result the hypotheses assume the multiplicativity of the functional $\phi$ at the cost of relaxing the linearity. Another key tool employed in the proof is a celebrated result in which linearity is not assumed in the hypotheses, we refer to the Kowlaski--S{\l}odkowski theorem.

\begin{theorem}\label{t KS}{\rm(Kowalski--S{\l}odkowski \cite{KoSlod})} 
Let $A$ be a complex Banach algebra, and let $\Delta : A\to \mathbb{C}$ be a functional satisfying $\Delta (0)=0$ and 
$$\Delta (a) - \Delta (b) \in \sigma (a-b),$$ for all $a,b\in A$. Then $\Delta$ is linear and multiplicative.
\end{theorem}

Up to the best of our knowledge, the Gleason--Kahana--\.{Z}elazko and the Kowalski--S{\l}odkowski theorems have not been properly treated in the case of Jordan--Banach algebras. Let us briefly justify their validity. 

\begin{theorem}\label{t GKZ JB} Let $\mathfrak{A}$ be a complex Jordan--Banach algebra, and let $\varphi : \mathfrak{A}\to \mathbb{C}$ be a non-zero linear functional. Then the following statements are equivalent:\begin{enumerate}[$(a)$]\item $\varphi (a)\in \hbox{J-}\sigma(a)$ for all $a\in \mathfrak{A}$.
\item $\varphi$ is unital or admits a unital extension to the unitization of $\mathfrak{A}$ and maps invertible elements in $\mathfrak{A}$ to non-zero complex numbers.
\item $\varphi$ is (Jordan-)multiplicative.
\end{enumerate}
\end{theorem}

\begin{proof} We shall only show that $(a)\Rightarrow (c)$, the other implications are more or less standard. There is no loss of generality in assuming that $\mathfrak{A}$ is unital. Let us fix an element $a\in \mathfrak{A}$ and denote by $\mathfrak{B}$ the Jordan--Banach algebra generated by $a$ and the unit element. It is well known that $\mathfrak{B}$ is a commutative unital complex Banach algebra. By the basic properties of the spectrum, the restricted functional $\varphi|_{\mathfrak{B}} : \mathfrak{B}\to \mathbb{C}$ satisfies $\varphi (a) \in \hbox{J-}\sigma_{\mathfrak{A}} (a) \subseteq \hbox{J-}\sigma_{\mathfrak{B}} (a)$, and in $\mathfrak{B}$ the Jordan and the usual spectrum of each element define the same set. It then follows from the Gleason--Kahane--\.{Z}elazko theorem for Banach algebras that  $\varphi|_{\mathfrak{B}}$ is multiplicative, and hence $\varphi(a^2) = \varphi(a)^2$ for all $a\in \mathfrak{A}$. It is routine to check that $\varphi$ is a (Jordan-)multiplicative functional. 
\end{proof}

The Jordan version of the Kowalski--S{\l}odkowski theorem cannot be obtained from the corresponding result for Banach algebras by restricting to the Jordan algebra, $\mathfrak{B},$ generated by two arbitrary elements and applying the Shirshov--Cohn theorem, essentially because this result is a purely algebraic and we do not know, in general, whether $\mathfrak{B}$ can be embedded as a Jordan subalgebra of a Banach algebra. However the original proof of Kowalski and S{\l}odkowski literary works in this case too (actually the spherical variants of the Kowalski--S{\l}odkowski theorem in \cite{LiPeWanWan19,Oi} also work for Jordan--Banach algebras, details are left to the reader). 
 
\begin{theorem}\label{t KS JB} Let $\mathfrak{A}$ be a complex Jordan--Banach algebra, and let $\Delta : \mathfrak{A}\to \mathbb{C}$ be a functional satisfying $\Delta (0)=0$ and $$\Delta (a) - \Delta (b) \in \hbox{J-}\sigma (a-b),$$ for all $a,b\in \mathfrak{A}$. Then $\Delta$ is linear and multiplicative.
\end{theorem}

One of the novelties in our approach here is that our arguments do not require any application of the Kowalski--S{\l}odkowski theorem but a weaker algebraic lemma also due to these authors. Namely, if $A$ is a complex Banach algebra every $\mathbb{R}$-linear spectral-valued functional $\varphi: A\to \mathbb{C}$ is $\mathbb{C}$-linear \cite[Lemma 2.1]{KoSlod}. If $\mathfrak{A}$ is a complex Jordan--Banach algebra and $\varphi : \mathfrak{A}\to \mathbb{C}$ is an $\mathbb{R}$-linear spectral-valued functional, by restricting to the Jordan--Banach subalgebra, $\mathfrak{B}$, generated by a single element $a\in \mathfrak{A}$, we can conclude, as in the proof of Theorem~\ref{t GKZ JB}, that $\varphi|_{\mathfrak{B}}$ (and thus $\varphi$) is  $\mathbb{C}$-linear. We state this conclusion.

\begin{lemma}\label{l KS r-linear spectral valued for JB} Let $\mathfrak{A}$ be a complex Jordan--Banach algebra, and let $\varphi : \mathfrak{A}\to \mathbb{C}$ be an $\mathbb{R}$-linear spectral-valued functional. Then $\varphi$ is  $\mathbb{C}$-linear.	
\end{lemma}

After reviewing the state-of-the-art we shall next exploit the consequences of the Lie--Trotter formulae for unital complex Jordan--Banach algebras.   

\begin{lemma}\label{sp1}
Let $X$ be a complex Banach space, let $A$ be a unital complex Banach algebra, and 
let $\phi\colon X\to A$ be a mapping satisfying:
\begin{enumerate}[$(a)$]
\item
$\phi(0)=\mathbf{1}$,
\item[(b)]
$\phi(x+y)=\phi(x)\, \phi(y),$ for all $x,y\in X$,
\item[(c)]
$\phi$ is continuous.
\end{enumerate}
Then there exists a unique $\mathbb{R}$-linear map $\psi\colon X\to A$ such that
\[
\phi(x)=e^{\psi(x)},\ \hbox{ for all } x\in X.
\]
\end{lemma}

\begin{proof} Fix an arbitrary $x\in X$. Then the mapping $\phi_x\colon\mathbb{R}\to A$ defined by
\[
\phi_x(t)=\phi(t x)\quad \forall t\in\mathbb{R},
\]
is a continuous one-parameter group on $A$, and therefore
there exists a unique $\psi(x)\in A$ such that
$\phi_x(t)=e^{t\psi(x)}$ for each $t\in\mathbb{R}$, i.e.,
\begin{equation}\label{e951}
\phi(t x)=e^{t\psi(x)}\quad\forall t\in\mathbb{R}
\end{equation} (cf. \cite[Theorem 1.1.31]{CabGarPalVol1}). In particular, $\phi(x)=e^{\psi(x)}$.
	
Our next goal is to show that $\psi$ is $\mathbb{R}$-linear.
Let $x,y\in X$ and $\alpha,\beta\in\mathbb{R}$.
For each $t\in\mathbb{R}$, we have
$$\begin{aligned} e^{t\psi(\alpha\,x+\beta\,y)}
	&=
	\phi\bigl(t(\alpha\,x+\beta\,y)\bigr)=
	\phi\bigl((t\alpha)x+(t\beta)y\bigr)=
	\phi\bigl((t\alpha)x\bigr)\phi\bigl((t\beta)y\bigr) \\ 
	&= e^{t\alpha\psi(x)}e^{t\beta\psi(y)}=
	e^{t(\alpha\psi(x)+\beta\psi(y))},
\end{aligned} $$ witnessing that $\psi(\alpha\,x+\beta\,y) = \alpha\psi(x)+\beta\psi(y)$.
\end{proof}

The Jordan versions of the Lie--Trotter formula play now their role. 

\begin{lemma}\label{sp2}
Let $\mathfrak{A}$ be a unital complex Jordan--Banach algebra, and let $f\colon \mathfrak{A} \to\mathbb{C}$ be a functional such that
\begin{enumerate}[$(a)$]
\item $f(\mathbf{1})= {1}$,
\item $f(x\circ y)=f(x) f(y),$ for all $x,y\in \mathfrak{A}$,
\item $f$ is continuous at $\mathbf{1}$.
\end{enumerate}
Then there exists a unique $\mathbb{R}$-linear functional $\psi\colon  \mathfrak{A}\to\mathbb{C}$ such that
\[
f(e^x)=e^{\psi(x)}, \hbox{ for all } x\in \mathfrak{A}.
\]
\end{lemma}

\begin{proof} We shall show that we can apply the previous Lemma~\ref{sp1} to the functional $\phi :\mathfrak{A}\to \mathbb{C}$ defined by $\phi(x): =f(e^x)$. First we observe that the continuity of $f$ at $\unit$ passes to each element of the form $e^x.$ Namely, if $(y_n)\to e^x$, the sequence $(e^{-x} \circ y_n)$ converges in norm to $\unit$, and hence $f(e^{-x}) f(y_n) = f (e^{-x} \circ y_n) \to f(\unit) =1$. It follows from the hypotheses that $f$ is continuous at $e^x$, which in turn proves the continuity of $\phi$. Therefore, it remains to show that  $\phi (x+ y) = \phi (x) \phi (y)$ for all $x,y\in \mathfrak{A}$. For this purpose we apply the last Lie--Trotter formula in Corollary~\ref{c Li-Trotter formulae Jordan-Banach algeras} and the hypotheses on $f$ to get $$\begin{aligned}
		\phi(x)\phi(y)=f(e^x)f(e^y) &=
		f\bigl((e^{x/n})^n\bigr)f\bigl((e^{y/n})^n\bigr) =\bigl(f(e^{x/n})\bigr)^n\bigl(f(e^{y/n})\bigr)^n\\		
		&=
		\bigl(f(e^{x/n})f(e^{y/n})\bigr)^n=
		\bigl(f(e^{x/n}\circ e^{y/n})\bigr)^n\\
		&= f\bigl((e^{x/n}\circ e^{y/n})^n\bigr)\to f(e^{x+y})=\phi(x+y).
	\end{aligned} $$
\end{proof}

We can next state our conclusion on (Jordan-)multiplicative frunctionals on a unital complex Jordan--Banach.

\begin{theorem}\label{t Maounche Jordan mutiplicative}
Let $\mathfrak{A}$ be a unital complex  Jordan--Banach algebra, and let $f\colon \mathfrak{A}\to\mathbb{C}$ be a (non-necessarily linear) spectral-valued, (Jordan-)multiplicative functional, that is, 
$$ f(x)\in\text{\rm J-}\sigma(x), \hbox{ and } f(x\circ y)=f(x)  f(y),\quad\forall x,y\in \mathfrak{A}.$$
Then
there exists a unique multiplicative and continuous (complex) linear functional $\psi\colon \mathfrak{A}\to\mathbb{C}$
such that
\[
f(e^x)=\psi(e^x),\hbox{ for all } x\in \mathfrak{A}.
\]
\end{theorem}

\begin{proof} On one hand we have $f(\mathbf{1})\in 
\text{J-}\sigma(\mathbf{1})=\{1\},$ and thus $f(\mathbf{1})=1$.

The continuity of $f$ at $\mathbf{1}$ can be obtained by a standard argument (cf. \cite[Théor\`{e}me]{Maouche1996}). For completeness,  we observe that $$f (\unit + x) \subseteq \text{J-}\sigma(\mathbf{1} + x) = 1+ \text{J-}\sigma(x),$$ which implies that $f(\unit)- f (\unit + x) \in \text{J-}\sigma(-x),$ and hence $|f(\unit)- f (\unit + x) | \leq \|x\|$. 

Lemma~\ref{sp2} proves the existence of a real linear functional $\psi\colon \mathfrak{A}\to\mathbb{C}$ satisfying $f(e^x) = e^{\psi(x)}$ for all $x\in \mathfrak{A}$. The rest of the proof is devoted to show that 
\[
\psi(x)\in \text{J-}\sigma(x),\hbox{ for all } x\in \mathfrak{A}.
\] Namely, by construction and the spectral mapping theorem, we have \begin{equation}\label{e1449}
\exp\bigl(\psi(x)\bigr)=
f(e^x)\in \text{J-}\sigma
 \bigl(e^x\bigr)=
\exp\bigl(\text{J-}\sigma(x)\bigr), \hbox{ for all } x\in \mathfrak{A}.
\end{equation}
We fix $x\in \mathfrak{A}$ and an arbitrary $t\in\mathbb{R}^+$ such that 
\[
\vert t\psi(x)\vert<1\hbox{ and }
\Vert tx\Vert<1.
\] Then $\psi(tx)=t\psi(x)\in\mathbb{D}$ and $\text{J-}\sigma (tx)\subset\mathbb{D},$ and since the complex exponential is injective on $\mathbb{D}$, we deduce from \eqref{e1449} that $\psi(tx)\in\text{J-}\sigma (tx).$
Consequently,
\[
\psi(x)=
t^{-1}\psi(tx)\in
t^{-1} \text{J-}\sigma (tx)=
\text{J-}\sigma(x).
\] We can therefore apply Lemma~\ref{l KS r-linear spectral valued for JB} to deduce that $\psi$ is complex linear, and consequently, by the Jordan version of the Gleason--Kahane--\.{Z}elazko theorem (see Theorem~\ref{t GKZ JB}), $\psi$ is multiplicative. 

We finally observe that, since $\psi$ is multiplicative, bounded and linear we have 
\[
f(e^x)=e^{\psi(x)}=\psi(e^x),\hbox{ for all } x\in \mathfrak{A}.
\]
\end{proof}

Let us note that contrary to the arguments in \cite{Maouche1996}, the previous proof does not require the application of the Kowalski-S{\l}odkowski theorem. 

Let $A$ denote a unital complex associative Banach algebra with unit $\textbf{1}$. In this case the set, $A^{-1},$ of invertible elements in $A$ is a multiplicative subgroup of $A$, which is not, in general connected. The  principal component of $A^{-1}$ (i.e. the connected component of $A^{-1}$ containing $\unit$) will be denoted by  $A_{\textbf{1}}^{-1}$. It is known that \begin{equation}\label{eq principal component invertible Banach} A_{\mathbf{1}}^{-1} = \left\{ e^{a_1}\cdots e^{a_n} : a_i\in A, \ n\in \mathbb{N} \right\},
\end{equation} and it is precisely the least subgroup of $A^{-1}$ containing $\exp(A)$ (cf. \cite[Propositions 8.6 and 8.7]{BonDunCNA73}).

It is worth to note that we can rediscover Maouche's theorem in \cite{Maouche1996} from our previous result. Namely, suppose $f: A\to \mathbb{C}$ is a (non-necessarily linear) multiplicative functional such that $f (x) \in \sigma(x),$ for all $x\in A$, where $A$ is a unital complex Banach algebra. Clearly $f$ is Jordan-multiplicative, and since $\sigma(x) = \hbox{J-}\sigma(x)$ for all $x\in A$, Theorem~\ref{t Maounche Jordan mutiplicative} assures the existence of a continuous and multiplicative functional $\psi : A\to \mathbb{C}$ satisfying $f(e^x) = \psi(e^x)$ for all $x\in A$.  Theorem 1 in \cite{Ze68} implies that $\psi$ is multiplicative. It follows from \eqref{eq principal component invertible Banach} that $f$ and $\psi$ agree on $A^{-1}_{\unit}$.

We shall keep a consistent notation in the Jordan case. The principal component of the invertible elements in a unital complex Jordan--Banach algebra $\mathfrak{A}$ will be denoted by $\mathfrak{A}^{-1}_{\unit}$. A deep result by Loos in \cite{Loos96} establishes that $$\mathfrak{A}_{\mathbf{1}}^{-1} = \left\{ U_{\exp(a_1)}\cdots U_{\exp(a_n)} (\textbf{1}) : a_i\in \mathfrak{A}, \ n\in \mathbb{N} \right\},$$ consequently, $\mathfrak{A}_{\mathbf{1}}^{-1}$ is open and closed in $\mathfrak{A}^{-1}$ and each one of its connected components is analytically arcwise connected \cite[Corollary]{Loos96}. Jordan products of invertible elements are not, in general, invertible so, an expression of the form in \eqref{eq principal component invertible Banach} with the associative product replaced with the Jordan product is hopeless. 

A subset $\mathfrak{M}$ of invertible elements in a unital complex Jordan--Banach algebra $\mathfrak{A}$ is called a \emph{quadratic subset} if $U_{a} (b) \in \mathfrak{M}$ for all $a,b\in \mathfrak{M}$. It is shown in \cite[Lemma 2.5]{Pe2021} that $\mathfrak{A}^{-1}_{\unit}$ characterises as the least quadratic subset of $\mathfrak{A}^{-1}$ containing $\exp(\mathfrak{A})$.  

The previous brief description motivates the hypotheses in the next result.

\begin{lemma}\label{l f preserves U} 
Let $\mathfrak{A}$ be a unital complex Jordan--Banach algebra, and let $f\colon \mathfrak{A}\to\mathbb{C}$ be a non-zero (non-necessarily linear) functional satisfying:\begin{enumerate}[$(a)$]
	\item $f$ is continuous at $\mathbf{1}$. 
	\item $f(U_x (y))=U_{f(x)} f(y),$ for all $x,y\in \mathfrak{A}$.
\end{enumerate}
Then there exists a unique real linear functional $\psi: \mathfrak{A}\to \mathbb{C}$ satisfying one of the following mutually exclusive statements \begin{enumerate}[$(1)$]
	\item $f(e^x)=e^{\psi(x)}$ for all $x\in \mathfrak{A}$. 
	\item $f(e^x)=- e^{\psi(x)}$ for all $x\in \mathfrak{A}$.
\end{enumerate} Furthermore, the first (respectively, the second) statement holds precisely when $f(\unit) = 1$  
(respectively, $f(\unit) = -1$).
\end{lemma}

\begin{proof} We define a functional $\phi : \mathfrak{A}\to \mathbb{C}$ by  $\phi(x)=f(e^x)$. We begin by showing some of the properties of $f$ and $\phi$. Since $$\phi (0) = f (e^0) = f (\unit) = U_{f (\unit)} (f (\unit)) = f (\unit)^3 =\phi(0)^3,$$ and $f$ is non-zero, it follows that $\phi(0)  =f(\unit)\in \{
	\pm 1\}$.
	
Assume first that $f(\unit) = 1$. In such a case, $f(x^2) = f\left(U_{x} (\unit) \right) = U_{f(x)} f(\unit) = f(x)^2$, and inductively, $f(x^n) = f(x)^n$ for all natural $n$ and $x\in \mathfrak{A}$. Furthermore, $\phi(x)\neq 0$ and $\phi(-x) = f(e^{-x}) = f(e^x)^{-1} = \phi (x)^{-1}$ for all $x\in \mathfrak{A}.$

If $(x_n)\to x$ in norm, the sequence $(U_{e^{\frac{x_n}{2}}} (e^{-x} ) )$ tends to $\unit$ in norm. By the hypotheses on $f$ we have $$\phi(x_n) \, \phi(x)^{-1}= f(e^{\frac{x_n}{2}} )^2 f(e^{-x})  = f(e^{{x_n}} ) f(e^{-x}) = f(e^{\frac{x_n}{2}} )^2 f(e^{-x}) = f\left( U_{e^{\frac{x_n}{2}}} (e^{-x}) \right) \to 1,$$ which proves the continuity of $\phi$.
	
Now we observe that \[
\bigl(U_{f(e^{x/n})}f(e^{y/n})\bigr)^n=f(e^{x/n})^{2n}f(e^{y/n})=f(e^{2x})f(e^y)=\phi(2x)\phi(y)
\]
and, on the other hand, by the second formula in Corollary~\ref{c Li-Trotter formulae Jordan-Banach algeras} combined with the continuity of $f$ at $e^{2x +y},$ we arrive to 
\[
\bigl(U_{f(e^{x/n})}f(e^{y/n})\bigr)^n=
\bigl(f\bigl(U_{e^{x/n}}e^{y/n}\bigr)\bigr)^n=
f\bigl(\bigl(U_{e^{x/n}}e^{y/n}\bigr)^n\bigr)\to f(e^{2x+y})=\phi(2x+y).
\]  Both identities together give  $\phi(2x+y) = \phi(2x)\phi(y),$ equivalently, $\phi(x+y) = \phi(x) \phi(y)$, for all $x,y\in \mathfrak{A}$. 

We are in a position to apply Lemma~\ref{sp1} to $\phi$ to deduce the existence of a (unique) real linear functional $\psi: \mathfrak{A}\to \mathbb{C}$ satisfying $f(e^x) = \phi (x) = e^{\psi(x)},$ for all $x\in \mathfrak{A}$. 

If $f(\unit) = -1$, the functional $-f$ satisfies the same hypotheses of $f$ and sends $\unit$ to $1$. By the conclusion in the first part, there exists a (unique) real linear functional $\psi: \mathfrak{A}\to \mathbb{C}$ satisfying $-f (x)  = e^{\psi(x)},$ for all $x\in \mathfrak{A}$.
\end{proof}

If $A$ and $B$ are unital Banach algebras and $\Delta: A\to B$ is a unital (non-necessarily linear) mapping, it is not hard to see that $\Delta$ is Jordan-multiplicative (i.e. $\Delta (a  \circ b) = \Delta (a)\circ \Delta(b)$) if and only if $\Delta \left(U_a (x)\right) = U_{\Delta \left(a\right)} \Delta \left(x\right)$ (just apply that, in this particular setting, $U_a (x) = a x a,$ for all $a,x\in A$).  If $\Delta: \mathfrak{A}\to \mathfrak{B}$ is a unital (non-necessarily linear) mapping we have 
$$ \Delta \left(U_a (x)\right) = U_{\Delta \left(a\right)} \Delta \left(x\right), \ \forall a,x\in\mathfrak{A} \Rightarrow \Delta \hbox{ is Jordan-multiplicative}.$$ The reverse implication is clearly true when $\Delta$ is linear, however, without linearity we do not know if the reverse implication holds. 

The special algebraic description of the principal component of the invertible elements in a unital complex Jordan--Banach algebra commented before \cite{Loos96}, together with the previous discussion, invite us to consider some special, but natural, hypotheses in our next result.  

\begin{theorem}\label{t Maouche Jordan Banach with natural hypotheses}
Let $\mathfrak{A}$ be a unital complex Jordan--Banach algebra, and let $f\colon \mathfrak{A}\to\mathbb{C}$ be a (non-necessarily linear) functional satisfying:\begin{enumerate}[$(a)$]
	\item $ f(U_x (y))=U_{f(x)}f(y),$ for all $x,y\in \mathfrak{A}$,
	\item $f(x)\in\text{\rm J-}\sigma (x),$  for all $x\in \mathfrak{A}.$
\end{enumerate}
Then there exists a unique continuous (Jordan-)multiplicative linear functional $\psi\colon \mathfrak{A}\to\mathbb{C}$
such that $ f(x)=\psi(x),$ for every $x$ in the connected component of set of all invertible elements of $\mathfrak{A}$ containing the unit element.
\end{theorem}

\begin{proof} As in the proof of Theorem~\ref{t Maounche Jordan mutiplicative}, $f(\unit) =1$, and the condition $1-f(\unit + x) \in \text{\rm J-}\sigma (x)$, implies that $|1-f(\unit + x)|\leq \|x\|$, which suffices to show that $f$ is continuous at $\unit$. By Lemma~\ref{l f preserves U}, there exists a real linear functional $\psi: \mathfrak{A}\to \mathbb{C}$ satisfying $f(e^x) = e^{\psi(x)}$ for all $x\in \mathfrak{A}$. Arguing as in the second part of the proof of Theorem~\ref{t Maounche Jordan mutiplicative}, we get $\psi(x) \in \text{\rm J-}\sigma (x)$. Lemma~\ref{l KS r-linear spectral valued for JB} proves that $\psi$ is $\mathbb{C}$-linear, and thus $f(e^x) = e^{\psi(x)} = \psi(e^x),$ for all $x\in \mathfrak{A}$. 
	
Finally, by the main result in \cite{Loos96}, the principal component of $\mathfrak{A}^{-1}$ is the set $$\mathfrak{A}_{\mathbf{1}}^{-1} = \left\{ U_{\exp(a_1)}\cdots U_{\exp(a_n)} (\textbf{1}) : a_i\in \mathfrak{A}, \ n\in \mathbb{N} \right\}.$$ Therefore, by combining the property  $f(e^x) = e^{\psi(x)} = \psi(e^x),$ for all $x\in \mathfrak{A}$ with hypothesis $(a)$, we conclude that $f(x) = \psi (x)$ for all $x\in \mathfrak{A}_{\mathbf{1}}^{-1}$.
\end{proof}

It should be noted that the proof of Theorem~\ref{t Maouche Jordan Banach with natural hypotheses} does not require any application of the Kowalski-S{\l}odkowski theorem but the weaker technical Lemma~\ref{l KS r-linear spectral valued for JB}.

\begin{remark}\label{remark 1-homogeneity} Let $\mathfrak{A}$ be a unital complex Jordan--Banach algebra, and let $f: \mathfrak{A}\to \mathbb{C}$ be a (non-necessarily linear) functional satisfying the hypotheses in Theorem~\ref{t Maouche Jordan Banach with natural hypotheses}. Then $f$ is $1$-homogeneous, that is,  $f(\lambda x)=\lambda f(x)$, for all $\lambda \in \mathbb{C}$. Indeed, 
	\[
	f(U_{\sqrt{\lambda} \mathbf{1}}(x))=U_{f(\sqrt{\lambda}\mathbf{1})}f(x)=
	f(\sqrt{\lambda}\mathbf{1})^2f(x)=\sqrt{\lambda}^2f(x)=\lambda f(x).
	\]
\end{remark}

\section{Continuous spectral-valued multiplicative maps}\label{sec: Brits}

We already know from \cite[Exemple]{Maouche1996} about the existence of a non-linear (and non-continuous) spectral-valued multiplicative functional $f : C[0,1]\to \mathbb{C}$. So, the conclusions in Theorems~\ref{t Maounche Jordan mutiplicative} and \ref{t Maouche Jordan Banach with natural hypotheses} are optimal. However, by adding the hypothesis of continuity, Problem~\ref{Problem reciprocal of the GKZ towards linearity} on the automatic linearity of spectral-valued multiplicative maps from a complex Banach algebra, has been recently solved by Brits, Mabrouk, Schulz, and Tour{\'e} for commutative unital C$^*$-algebras, von Neumann algebras, and for unital C$^*$-algebras (cf. \cite{TouShulzBrits2017, TouBrits2020,BritsMabroukToure2021}).  One of the key tools employed in the known solutions for unital C$^*$-algebras is a tool, developed in \cite[Lemma 3.6]{TouShulzBrits2017}, showing that for each continuous multiplicative spectral-valued functional $f$ from a unital C$^*$-algebra $A$ into $\mathbb{C}$  and each  $x \in A$, the function $\lambda\to f(x- \lambda  \unit)$ is analytic on the unbounded component of $\mathbb{C}\backslash \sigma(x)$. In our next result we establish an extension of the just quoted tool by showing that, even under more general hypotheses, the commented function is in fact affine.

\begin{lemma}\label{sp20}
	Let $\mathfrak{A}$ be a unital Jordan--Banach algebra, and let $f\colon \mathfrak{A}\to\mathbb{C}$ be a (non-necessarily linear) functional satisfying: 
	\begin{enumerate}[$(a)$]
		\item $f$ is continuous;
		\item $f(U_x (y))=U_{f(x)}f(y)$ for all $x,y\in \mathfrak{A}$;
		\item $f(x)\in\text{\rm J-}\sigma(x)$ for each $x\in \mathfrak{A}$.
	\end{enumerate}
	Let $\psi$ be the bounded (Jordan-)multiplicative linear functional given by Theorem~\ref{t Maouche Jordan Banach with natural hypotheses}.
	If $x\in \mathfrak{A}$ and $\mathcal{O}_x$ denotes the unbounded component of 
	$\mathbb{C}\setminus\text{\rm J-}\sigma(x)$, then
	\begin{equation*}
		f(\lambda\mathbf{1}-x)=\lambda-\psi(x),\hbox{ for all }\lambda\in \mathcal{O}_x.
	\end{equation*}
	If $x\in J$ is such that
	$\text{\rm J-}\sigma(x)$ is contained in a line, then
	\begin{equation*}
		f(\lambda\mathbf{1}-x)=\lambda-\psi(x),\hbox{ for all } \lambda\in\mathbb{C},
	\end{equation*}
	and, in particular, 
	\[
	f(x)=\psi(x).
	\]
\end{lemma}

\begin{proof} We observe that 
	\begin{equation}\label{e1813}
		\lambda\mathbf{1}-x\in \mathfrak{A}^{-1}_{\unit},\quad\forall \lambda\in \mathcal{O}_x\setminus\{0\}.
	\end{equation} The desired statement can be directly deduced from \cite[Remark 2.4$(2)$]{Pe2021}, however, we include a short argument here for completeness reasons.  	We observe that $\lambda\mapsto \mathbf{1}-\lambda^{-1}x$ is a continuous mapping from the
	domain $\mathcal{O}_x\setminus\{0\}$ to $\mathfrak{A}^{-1},$ and hence the set
	\[
	\bigl\{\mathbf{1}-\lambda^{-1}x\colon \lambda\in \mathcal{O}_x\setminus\{0\}\bigr\}\subseteq \mathfrak{A}^{-1} 
	\]
	is connected. On the other hand, $\mathbf{1}$ lies in the closure of the above set as
	\[
	\lim_{\substack{\vert\lambda\vert\to\infty\\ \lambda\in \mathcal{O}_x\setminus\{0\}}}\mathbf{1}-\lambda^{-1}x=\mathbf{1}.
	\]
	Thus,
	\[
	\bigl\{\mathbf{1}-\lambda^{-1}x\colon \lambda\in \mathcal{O}_x\setminus\{0\}\bigr\} \ \cup \ \bigl\{\mathbf{1}\bigr\} \subseteq \mathfrak{A}^{-1} 
	\]
	is connected and contains $\mathbf{1}$. Hence the above set is contained in $\mathfrak{A}_{\unit}^{-1}$. 
	By taking into account that $\mathbb{C}^* \, \mathfrak{A}_{\unit}^{-1}  \subseteq \mathfrak{A}_{\unit}^{-1},$ it may be concluded that \eqref{e1813} holds.
	Theorem~\ref{t Maouche Jordan Banach with natural hypotheses} assures that 
	\[
	f(\lambda\mathbf{1}-x)=\psi(\lambda\mathbf{1}-x)=\lambda\psi(\mathbf{1})-\psi(x)=\lambda-\psi(x),
	\quad\forall \lambda\in \mathcal{O}_x\setminus\{0\}.
	\]
	In the case where $0\in\mathcal{O}_x$, the continuity of $f$ allows to extend the above identity to 
	the whole domain $\mathcal{O}_x$.
	
	Suppose that $\text{\rm J-}\sigma(x)$ is contained in a line. Then $\mathcal{O}_x$ is dense in $\mathbb{C}$
	and the required indentity follows from the previously proved identity together with the continuity of $f$.
	By taking $\lambda=0$, we get $f(x)=\psi(x)$.
\end{proof}

The following corollary, which is a strengthened version of \cite[Lemma 3.6]{TouShulzBrits2017} with more general hypotheses and a stronger conclusion, follows from the previous lemma together with the commented result by \.{Z}elazko asserting that each Jordan-multiplicative linear functional on a Banach algebra is multiplicative \cite{Ze68}. 

\begin{corollary}\label{c generalization of TSB2017 L3.6}
Let ${A}$ be a unital Banach algebra, and let $f\colon {A}\to\mathbb{C}$ be a (non-necessarily linear) functional satisfying: 
\begin{enumerate}[$(a)$]
	\item $f$ is continuous;
	\item $f$ is multiplicative;
	\item $f(x)\in \sigma(x)$ for each $x\in {A}$.
\end{enumerate}
Let $\psi$ be the bounded (Jordan-)multiplicative linear functional given by Maouche's theorem (see Theorem~\ref{t Maouche}).
If $x\in {A}$ and $\mathcal{O}_x$ denotes the unbounded component of 
$\mathbb{C}\setminus\sigma(x)$, then
\begin{equation*}
	f(\lambda\mathbf{1}-x)=\lambda-\psi(x),\hbox{ for all }\lambda\in \mathcal{O}_x.
\end{equation*} In particular, the mapping $\lambda\mapsto f(\lambda\mathbf{1}-x)$ is affine on  $\mathcal{O}_x$.
If $x\in J$ is such that $\sigma(x)$ is contained in a line, then
\begin{equation*}
	f(\lambda\mathbf{1}-x)=\lambda-\psi(x),\hbox{ for all } \lambda\in\mathbb{C},
\end{equation*}
and, in particular,  $f(x)=\psi(x).$	
\end{corollary}

We can now turn our attention to the second question in Problem~\ref{Problem reciprocal of the GKZ towards linearity Jordan-Banach} in the case of unital JB$^*$-algebras. According to the standard notation, henceforth we shall write $\mathfrak{A}_{sa}$ for the set of all self-adjoint elements in $\mathfrak{A}$.  It is known that the JB$^*$-subalgebra of $\mathfrak{A}$ generated by a single element $a\in \mathfrak{A}_{sa}$ is isometrically isomorphic as JB$^*$-algebra to a commutative C$^*$-algebra (cf. \cite[Theorem 3.2.2]{HOS}).

\begin{theorem}\label{t BMT for JB}
	Let $\mathfrak{A}$ be a unital JB$^*$-algebra, and let $f\colon \mathfrak{A} \to\mathbb{C}$ be a (non-necessarily linear) functional satisfying:
	\begin{enumerate}[$(a)$]
		\item $f$ is continuous;
		\item $f(U_x (y))=U_{f(x)}f(y),$ for all $x,y\in \mathfrak{A}$;
		\item $f(x)\in\text{\rm J-}\sigma(x),$ for each $x\in \mathfrak{A}$.
	\end{enumerate}
	Then $f$ is a linear (Jordan-)multiplicative functional on $\mathfrak{A}$.
\end{theorem}

\begin{proof} Let $\psi\in \mathfrak{A}^*$ be the functional given  by Theorem~\ref{t Maouche Jordan Banach with natural hypotheses}.  We begin by observing that
	\begin{equation}\label{e1955}
		f(x+\lambda\mathbf{1})=\psi(x)+\lambda,\quad\forall x\in \mathfrak{A}_{sa}, \ \lambda\in \mathbb{C},
	\end{equation}
	and particularly
	\begin{equation}\label{e1956}
		f(x)=\psi(x),\quad\forall x\in \mathfrak{A}_{sa}.
	\end{equation}
	Indeed, if $x\in \mathfrak{A}_{sa}$, then $\text{J-}\sigma(x)\subset\mathbb{R}$ and 
	Lemma~\ref{sp20} gives the claimed identities.
	
Now, let $x,y\in \mathfrak{A}_{sa}$ and assume that $x$ is invertible and positive. Then we have
	\begin{equation}\label{e2041}
		f(x+i\, y)=\psi(x+i\,y).
	\end{equation}
To see this, we observe that $U_{x^{-1/2}}(y)\in \mathfrak{A}_{sa},$ and thus the hypotheses and \eqref{e1956} 
give $$\begin{aligned}
f(x+i\,y) &= f\bigl(U_{x^{1/2}}(\mathbf{1}+i\,U_{x^{-1/2}}(y))\bigr)= U_{f(x^{1/2})} f\bigl(\mathbf{1}+i\,U_{x^{-1/2}}(y)\bigr)\\&=
	U_{\psi(x^{1/2})}\psi\bigl(\mathbf{1}+i\,U_{x^{-1/2}}(y)\bigr)=
	\psi(x)\bigl(1+\psi(i\,U_{x^{-1/2}}(y))\bigr)\\
	&=
	\psi(x)\bigl( 1 + i\,U_{\psi(x^{-1/2})}\psi(y)\bigr)=
	\psi(x)\bigl({1}+i\,\psi(x^{-1})\psi(y)\bigr)\\
	&=
	\psi(x)+i\,\psi(y)=\psi(x+i\,y).
\end{aligned}$$

Let $x,y\in \mathfrak{A}_{sa}$ and assume merely that $x\ge 0$. We shall show that
	\begin{equation}\label{e2042}
		f(x+i\, y)=\psi(x+i\,y).
	\end{equation}
By applying \eqref{e2041} with $x$ replaced by $n^{-1}\mathbf{1}+x$ we obtain
	\[
	f(n^{-1}\mathbf{1}+x+i\,y)=\psi(n^{-1}\mathbf{1}+x+i\,y)
	\]
and then by taking limits in $n$ and having in mind the continuity of both $f$ and $\psi$ we get \eqref{e2042}.
	
We pick now $x,y\in \mathfrak{A}_{sa}$, with $f(x)\ne 0$. We shall prove that $f(x+i\,y)=\psi(x+i\,y).$ Write 
	$x=a-b$ with $a,b\in \mathfrak{A}_{sa}$, $a,b\ge 0,$ and $U_a(b)= U_b (a) =0$ (this can be done by just observing that the JB$^*$-subalgebra of $\mathfrak{A}$ generated by $x$ is JB$^*$-isomorphic to a commutative C$^*$-algebra \cite[Theorem 3.2.2]{HOS}). Then,
	taking into account the hypotheses,  \eqref{e2042}, and the fact that $U_a (a)\ge 0$, we arrive to the following identity
$$\begin{aligned}
		f(a)^2 f(x+i\,y) &=U_{f(a)} f(x+i\,y)=
	f\bigl(U_a(x+i\,y)\bigr)=
	f\bigl(U_a (x)+i\,U_a (y)\bigr) \\
	&=	f\bigl(U_a (a) + i\,U_a (y))=
	\psi(U_a (a)+i\,U_a (y))\\
	&=
	\psi(U_a (x) +i\,U_a (y))=
	\psi(U_a(x+i\,y))\\
	&=
	U_{\psi(a)}\psi(x+i\,y)=
	\psi(a)^2\psi(x+i\,y)=
	f(a)^2 \psi(x+i\,y).
\end{aligned}
$$
	Similarly, we obtain
$$\begin{aligned}
f(b)^2f(x+i\,y) &= U_{f(b)}f(x+i\,y)=
	f\bigl(U_b(x+i\,y)\bigr)=
	f\bigl(U_b (x)+i\,U_b (y)\bigr)\\ 
	&=	f\bigl(-U_b (b) +i\,U_b (y))=
	-f\bigl(U_b (b)-i\,U_b (y))\\
	&=
	-\psi(U_b (b)-i\,U_b (y)) =	\psi(-U_b (b)+i\, U_b (y)) \\
	&=
	\psi(U_b (x) +i\,U_b (y))=
	\psi(U_b(x+i\,y))=U_{\psi(b)}\psi(x+i\,y) \\
	&=
	\psi(b)^2\psi(x+i\,y)=
	f(b)^2\psi(x+i\,y).
\end{aligned} 
$$
We thus get
	\begin{equation}\label{eq fa2 + fb2 times} \bigl[
		f(a)^2+f(b)^2
		\bigr]
		\bigl[
		f(x+i\,y)-\psi(x+i\,y)
		\bigr]=0.
	\end{equation}
Now we claim that $f(a)^2+f(b)^2\neq 0$. Otherwise, since $f(a)=\psi(a),f(b)=\psi(b)\in\mathbb{R}$ we derive that $f(a)=f(b)=0$, and therefore 
\[f(x)=\psi(x)=\psi(a)-\psi(b)=f(a)-f(b)=0,
\] contradicting that $f(x) \neq 0$.  It follows from \eqref{eq fa2 + fb2 times} that 
	\[
	f(x+i\,y)=\psi(x+i\,y),
	\] 
	as required.

Finally, assume that $f(x)=0$. Then, for each $n\in\mathbb{N}$, 
	$f(x+n^{-1} \mathbf{1})=\psi(x+n^{-1} \mathbf{1}) =\psi(x) + n^{-1} \psi( \mathbf{1}) =f(x) + n^{-1} =n^{-1}\ne 0$, and by applying what has been previously proved
	we obtain
	\[
	f(x+n^{-1}\mathbf{1}+i\,y)=\psi(x+n^{-1}\mathbf{1}+i\,y).
	\]
	By taking limits and using the continuity of both $f$ and $\psi$ we obtain that these two maps coincide on the whole $\mathfrak{A}$.
\end{proof}

Clearly, we can rediscover the main result in \cite{BritsMabroukToure2021} by combining the previous theorem with \.{Z}elazko's theorem in \cite{Ze68}.

\begin{corollary}{\rm\cite[Theorem 2.1]{BritsMabroukToure2021}} Let $A$ be a unital C$^*$-algebra. Then any continuous multiplicative functional $f:A\to \mathbb{C}$ satisfying $f(a)\in \sigma (a),$ for each $a\in A,$ is linear and hence a character of $A$.
\end{corollary}

It remains as an open problem whether the conclusions in Theorem~\ref{t Maouche}, Theorem~\ref{t Maouche Jordan Banach with natural hypotheses}, \cite[Theorem 2.1]{BritsMabroukToure2021}, and  Theorem~\ref{t BMT for JB} remain true in the non-unital case. 

\medskip\medskip

\textbf{Acknowledgements}\medskip

Second and third authors partially supported by
MCIN/AEI/10.13039/501100011033 and by ``ERDF A way of making Europe'' grant PID2021-122126NB-C31
and Junta de Andalucía grants PY20$\underline{\ }$00255, FQM185, and FQM375.
First author supported by grant FPU21/00617 at University of Granada founded by Ministerio de Universidades (Spain), 
and by the IMAG--Mar{\'i}a de Maeztu grant CEX2020-001105-M/AEI/10.13039/ 501100 011033

\end{document}